 \documentclass[twoside,reqno,11pt]{fcaa-var}

\usepackage{graphicx}
\usepackage{epsfig}
\usepackage{amsthm}
\usepackage{amsmath}
\usepackage{latexsym}
\usepackage{amsfonts}
\usepackage{amssymb}

\newtheoremstyle{theorem}
  {15pt}          
  {15pt}  
  {\sl}  
  {\parindent}
  {\sc}  
  {. }   
  { }    
  {}     
\theoremstyle{theorem}

\newtheoremstyle{defi}
  {15pt}          
  {15pt}  
  {\rm}  
  {\parindent}     
  {\sc}  
  {. }    
  { }    
  {}     
\theoremstyle{defi}

 \def\proofname{\indent\rm P\ r\ o\ o\ f.}
 \def\proofend{\hfill$\Box$}
 \def\qed{\hfill$\Box$} 

 \usepackage{upref}

 \def\theequation{\arabic{section}.\arabic{equation}}

 \newtheorem{theo}{Theorem}[section]
 \newtheorem{lem}{Lemma}[section]
 \newtheorem{prop}{Proposition}[section]
 \newtheorem{cor}{Corollary}[section]
 \newtheorem{defi}{Definition}[section]

 \def\Lip{\operatorname{Lip}}
 \def\C{C}
 \def\Lp{L^p}
 \def\Linf{L^\infty}
 \def\diam{\operatorname{diam}}
 \newcommand{\argmax}[1]{\underset{#1}{\operatorname{argmax} \,}}
 \newcommand{\argmin}[1]{\underset{#1}{\operatorname{argmin} \,}}
 \newcommand{\esssup}[1]{\underset{#1}{\operatorname{ess\,sup} \,}}

 \def\themycopyrightfootnote{\vspace*{3pt}
 \copyright \, Year\,  Diogenes  Co., Sofia
   \par  \noindent pp. xxx--xxx, DOI: ..................
   \hfill  \vspace*{-36pt}
  }

 \title[FRACTIONAL DERIVATIVES OF LYAPUNOV FUNCTIONS]
 {FRACTIONAL DERIVATIVES OF CONVEX LYAPUNOV FUNCTIONS AND CONTROL PROBLEMS IN FRACTIONAL ORDER SYSTEMS \\ [4pt] IN ``FCAA'' JOURNAL}
 \author[M.I. Gomoyunov]{Mikhail I. Gomoyunov $^1$}

\begin{document}

 \vbox to 2.5cm { \vfill }


 \bigskip \medskip

 \begin{abstract}

The paper is devoted to the development of control procedures with a guide for conflict-controlled dynamical systems described by ordinary fractional differential equations with the Caputo derivative of an order $\alpha \in (0, 1).$ For the case when the guide is in a certain sense a copy of the system, a mutual aiming procedure between the initial system and the guide is elaborated. The proof of proximity between motions of the systems is based on the estimate of the fractional derivative of the superposition of a convex Lyapunov function and a function represented by the fractional integral of an essentially bounded measurable function. This estimate can be considered as a generalization of the known estimates of such type. An example is considered which illustrates the workability of the proposed control procedures.

 \medskip


 \smallskip


 \end{abstract}

 \maketitle

 \vspace*{-16pt}


\section*{Introduction}

Control procedures 
with auxiliary guides play an important role in the positional differential games theory (see, e.g., \cite{NNKrasovskii_ANKrasovskii_1995, NNKrasovskii_AISubbotin_1988}). They are used both in obtaining theoretical results and in developing numerically realizable and robust optimal control schemes for different types of dynamical systems, including systems with disturbances and counteractions (see, e.g., \cite{YuVAverboukh_2017,  NNKrasovskii_ANKotelnikova_2012, NYuLukoyanov_ARPlaksin_2015, VIMaksimov_2016, ARMatviychuk_VNUshakov_2006}). The goal of this paper is to develop such control procedures with a guide for fractional order conflict-controlled systems which motion is described by ordinary differential equations with the Caputo fractional derivative of an order $\alpha \in (0, 1).$ In the paper, a case is considered when a motion of the guide is described by the similar differential equations.

One of the main difficulties in design of control procedures with a guide is to ensure the proximity between motions of the initial system and the guide. The most useful tool here is the Lyapunov functions technique. When trying to extend the results obtained for the first order systems to the fractional order ones, a well-known problem arises that involves calculating the fractional derivative of the superposition of a Lyapunov function and a system motion. In \cite{NAguila-Camacho_MADuarte-Mermoud_JAGallegos_2014, AAAlikhanov_2010} the upper bound for such derivative was obtained for a quadratic Lyapunov function. Later, similar inequalities were proved for more general classes of convex Lyapunov functions (see, e.g., \cite{WChen_HDai_YSong_ZZhang_2017} and references therein). However, the validity of these estimates was established under certain assumptions about smoothness of a motion (at least, absolutely continuity). Moreover, these differentiability properties were essentially used in the proofs. Thus, a system motion is required to be smooth enough in order to these estimates can be applied.

On the other hand, for the considered in the paper conflict-controlled systems, it is natural that the right-hand side of the closed-loop system depends on the time variable explicitly. This leads to the fact that a system motion does not have to be differentiable. Indeed, there exist (see, e.g., \cite{BRoss_SGSamko_ERLove_1994}) nowhere differentiable functions that have continuous fractional derivatives of any order $\alpha \in (0, 1).$ Consequently, one can consider these functions as the solutions of the simplest fractional order equations with the continuous right-hand side that depends only on the time variable. However, for these solutions, the estimates from \cite{NAguila-Camacho_MADuarte-Mermoud_JAGallegos_2014, AAAlikhanov_2010, WChen_HDai_YSong_ZZhang_2017} can not be used.

For the controlled fractional order systems, in \cite{DIdczak_RKamocki_2011} it was proposed to consider a motion of the system as a function represented by the fractional order Riemann-Liouville (R.-L.) integral of a summable function, without requiring any differentiability properties, and the existence and uniqueness of such a motion were proved. In the paper this notion of a motion is used. But, due to the stronger assumptions on the right-hand side of the motion equation, this notion is slightly modified: instead of summable functions, measurable essentially bounded functions are considered. The corresponding existence and uniqueness results are given in Theorem~\ref{Thm_existence}. In order to apply discussed above estimates, concerning Lyapunov functions technique, for such motions, it was necessary to prove that the estimates from \cite{WChen_HDai_YSong_ZZhang_2017} are valid for functions represented by the R.-L. integral of measurable essentially bounded functions. This result is given in Lemma~\ref{Lem_differentiating_of_convex}.

This lemma constitutes the basis of the proof of proximity between a motion of the conflict-controlled fractional order system and the guide when a suitable mutual aiming procedure is used. It seems that these results may be applied for the development of the theory and numerical methods in positional differential games for fractional order systems. It should be noted also that some kinds of pursuit-evasion differential games in fractional order systems were considered earlier (see, e.g., \cite{AAChikrii_IIMatichin_2010, NNPetrov_2017}).

The paper is organized as follows. In Sect.~1 the definitions and some basic properties of the fractional order R.-L. integral, R.-L. and Caputo derivatives are given. Sect.~2 deals with a Cauchy problem for an ordinary differential equation with the Caputo fractional derivative of an order $\alpha \in (0, 1).$ The notion of a solution of this Cauchy problem is proposed, the existence and uniqueness of such a solution are proved. In Sect.~3 the estimate of the R.-L. fractional derivative of the superposition of a convex Lyapunov function and the solution of the Cauchy problem is obtained. The case when this solution is smooth (Lipshitz continuous) and the general case are considered separately. Sect.~4 deals with a conflict-controlled fractional order dynamical system. Basic notions and system motion properties are given, an auxiliary guide is introduced. In Sect.~5 the mutual aiming procedure that ensures proximity between motions of the initial system and the guide is proposed. The obtained results are illustrated by numerical simulations in Sect.~6. Concluding remarks are given in Sect~7.

\section{Notations, Definitions and Preliminary Results}

\setcounter{section}{1}
\setcounter{equation}{0}\setcounter{theorem}{0}

Let $k \in \mathbb{N}$ and $\mathbb{R}^k$ be the $k$-dimensional Euclidian space with the scalar product $\langle \cdot, \cdot \rangle$ and the norm $\|\cdot\|.$ By $B(r) \subset \mathbb{R}^k,$ $r \geqslant 0,$ we denote the closed ball with the center in the origin and the radius $r.$ Let $T > 0$ and the segment $[0, T] \subset \mathbb{R}$ be endowed with the Lebesgue measure. For $p \in [1, \infty),$ by $\Lp([0, T], \mathbb{R}^k),$ we denote the Banach space of (classes of equivalence of) $p$-th power integrable functions $x: [0, T] \rightarrow \mathbb{R}^k$ with the norm
\begin{equation*}
    \|x(\cdot)\|_p = \Big( \int_{0}^{T} \|x(t)\|^p dt \Big)^{1/p}.
\end{equation*}
By $\Linf([0, T], \mathbb{R}^k),$ we denote the Banach space of (classes of equivalence of) essentially bounded measurable functions $x: [0, T] \rightarrow \mathbb{R}^k$ with the norm
\begin{equation*}
    \|x(\cdot)\|_\infty = \esssup{t \in [0, T]} \|x(t)\|.
\end{equation*}
Let $\C([0, T], \mathbb{R}^k)$ be the Banach space of continuous functions $x: [0, T] \rightarrow~\mathbb{R}^k$ \linebreak with the norm $\|\cdot\|_\infty,$ $\Lip([0, T], \mathbb{R}^k) \subset \C([0, T], \mathbb{R}^k)$ be the set of Lipschitz continuous functions, $\Lip^0([0, T], \mathbb{R}^k) = \{x(\cdot) \in \Lip([0, T], \mathbb{R}^k): x(0) = 0\}.$

\subsection{Riemann-Liouville Fractional Order Integral}

\begin{defi}[see {\cite[Definition~2.1]{SGSamko_AAKilbas_OIMarichev_1993}}]
    For a function $\varphi: [0, T] \rightarrow~\mathbb{R}^k,$ the (left-sided) R.-L. fractional integral of an order $\alpha \in (0, 1)$ is defined by
    \begin{equation*}
        (I^\alpha \varphi) (t) = \frac{1}{\Gamma(\alpha)} \int_{0}^{t} \frac{\varphi(\tau)}{(t-\tau)^{1-\alpha}} d\tau, \quad t \in [0, T],
    \end{equation*}
    where $\Gamma(\cdot)$ is the Euler gamma-function (see, e.g., \cite[(1.54)]{SGSamko_AAKilbas_OIMarichev_1993}).
\end{defi}

Let us describe some properties of the R.-L. fractional integral.
\begin{prop} \label{Prop_integral_properties}
    Let $\alpha \in (0, 1)$ and $p \in (1/\alpha, \infty].$ Then:
    \begin{itemize}
        \item[($A.1$)] For any $\varphi(\cdot) \in \Lp ([0, T], \mathbb{R}^k),$ the value $(I^\alpha \varphi)(t)$ is well defined for any $t \in [0, T],$ and $(I^\alpha \varphi)(0) = 0.$
        \item[($A.2$)] There exists $H_p > 0$ such that, for any $\varphi(\cdot) \in \Lp ([0, T], \mathbb{R}^k)$ and any $t, \tau \in [0, T],$ the inequality below is valid:
            \begin{equation*}
                \|(I^\alpha \varphi) (t) - (I^\alpha \varphi) (\tau)\|
                \leqslant H_p \|\varphi(\cdot)\|_p |t - \tau|^{\alpha - 1/p},
            \end{equation*}
            where $1/p = 0$ if $p = \infty.$ In particular, $(I^\alpha \varphi)(\cdot) \in \C([0, T], \mathbb{R}^k)$ for any $\varphi(\cdot) \in \Lp([0, T], \mathbb{R}^k).$
        \item[($A.3$)] The operator $I^\alpha: \Lp ([0, T], \mathbb{R}^k) \rightarrow \C([0, T], \mathbb{R}^k)$ is linear and compact (i.e., maps bounded sets from $\Lp ([0, T], \mathbb{R}^k)$ into relatively compact sets from $\C([0, T], \mathbb{R}^k)$), and, in particular, is continuous.
        \item[($A.4$)] If $\varphi(\cdot) \in \Lip^0([0, T], \mathbb{R}^k),$ then $(I^\alpha \varphi) (\cdot) \in \Lip^0([0, T], \mathbb{R}^k).$
    \end{itemize}
\end{prop}
\proof
    Statements ($A.1$) and ($A.2$) are proved in \cite[Theorem~3.6, Remark~3.3]{SGSamko_AAKilbas_OIMarichev_1993} (see also \cite[Theorem~2.6]{KDiethelm_2010}). The validity of property ($A.3$) follows from ($A.2$) and Arcel\`{a}-Ascoli theorem (see, e.g., \cite[Ch.~I, \S~5, Theorem~4]{LVKantorovich_GPAkilov_1982}). The proof of statement ($A.4$) can be found in \cite[Theorem~3.1]{SGSamko_AAKilbas_OIMarichev_1993} (see also \cite[Theorem~2.5]{KDiethelm_2010}).
\proofend

Let us recall a fractional version of Bellman-Gronwall lemma.
\begin{lem}[see {\cite[Lemma~6.19]{KDiethelm_2010}}] \label{Lem_Bellman_Gronwall}
  Let $\varepsilon \geqslant 0,$ $\lambda \geqslant 0$ and a function $x(\cdot) \in \C([0, T], \mathbb{R})$ satisfy the inequality
  \begin{equation*}
        |x(t)|
        \leqslant \varepsilon + \frac{\lambda}{\Gamma(\alpha)} \int_{0}^{t} \frac{|x(\tau)|}{(t - \tau)^{1 - \alpha}} d\tau,
        \quad t \in [0, T].
  \end{equation*}
  Then the following inequalities hold:
  \begin{equation*}
    |x(t)| \leqslant \varepsilon E_\alpha (\lambda t^\alpha)
    \leqslant \varepsilon E_\alpha (\lambda T^\alpha), \quad t \in [0, T],
  \end{equation*}
  where $E_\alpha(\cdot)$ is the Mittag-Leffler function (see, e.g., \cite[(1.90)]{SGSamko_AAKilbas_OIMarichev_1993}).
\end{lem}

\subsection{Riemann-Liouville and Caputo Fractional Order Derivatives}

\begin{defi}[see {\cite[Definition~2.2]{SGSamko_AAKilbas_OIMarichev_1993}}]
    For a function $x: [0, T] \rightarrow~\mathbb{R}^k,$ the (left-sided) R.-L. fractional derivative of an order $\alpha \in (0, 1)$ is defined by
    \begin{equation*}
        (D^\alpha x) (t)
        = \frac{1}{\Gamma(1 - \alpha)} \frac{d}{dt} \int_{0}^{t} \frac{x(\tau)}{(t - \tau)^{\alpha}} d\tau, \quad t \in [0, T].
    \end{equation*}
\end{defi}

\begin{defi}[see {\cite[Definition~2.3]{SGSamko_AAKilbas_OIMarichev_1993}}]
    Let $I^\alpha (\Linf([0, T], \mathbb{R}^k))$ denote the set of functions $x: [0, T] \rightarrow \mathbb{R}^k$ represented by the R.-L. fractional integral of an order $\alpha \in (0, 1)$ of a function $\varphi(\cdot) \in \Linf([0, T], \mathbb{R}^k):$ $x(t) = (I^\alpha \varphi)(t),$ $t \in [0, T].$
\end{defi}

Let us describe some properties of the R.-L. fractional derivative.
\begin{prop} \label{Prop_differentiating_I^alpha}
    Let $\alpha \in (0, 1).$ If $x(\cdot) \in I^\alpha(\Linf ([0, T], \mathbb{R}^k)),$ then:
    \begin{itemize}
      \item[($B.1$)] The value $(D^\alpha x)(t)$ is well defined for almost every $t \in [0, T],$ and $(D^\alpha x)(\cdot) \in \Linf([0, T], \mathbb{R}^k).$
      \item[($B.2$)] The equality $\big( I^\alpha ( D^\alpha x) \big) (t) = x (t)$ is valid for any $t \in [0, T].$
      \item[($B.3$)] Let $\varphi(\cdot) \in \Linf([0, T], \mathbb{R}^k)$ be such that $x(t) = (I^\alpha \varphi)(t),$ $t \in [0, T].$ Then $\varphi(t) = (D^\alpha x)(t)$ for almost every $t \in [0, T].$
    \end{itemize}

    Moreover, if $x(\cdot) \in \Lip^0([0, T], \mathbb{R}^k),$ then:
    \begin{itemize}
      \item[($B.4$)] The value $(D^\alpha x) (t)$ is well defined for any $t \in [0, T],$ and the representation formula below holds:
            \begin{equation} \label{representation_formula}
                (D^\alpha x) (t) = \frac{1}{\Gamma(1 - \alpha)} \int_{0}^{t} \frac{\dot{x}(\tau)}{(t - \tau)^\alpha} d\tau,
                \quad t \in [0, T].
            \end{equation}
            where $\dot{x}(t) = d x(t)/dt,$ $t \in [0, T].$
      \item[($B.5$)] The following inclusions are valid: $x(\cdot) \in I^\alpha(\Linf([0, T], \mathbb{R}^k))$ and $(D^\alpha x)(\cdot) \in I^{1 - \alpha}(\Linf([0, T], \mathbb{R}^k)).$ In particular, $(D^\alpha x)(0) = 0.$
    \end{itemize}
\end{prop}
\proof
    Statements ($B.1$) and ($B.2$) are proved by the scheme from \cite[Theorem~2.4]{SGSamko_AAKilbas_OIMarichev_1993} (see also \cite[Theorem~2.22]{KDiethelm_2010}). The validity of property ($B.3$) follows from \cite[Lemma~2.1]{DIdczak_RKamocki_2011}. Statement ($B.4$) can be established by the scheme from \cite[Lemmas~2.1, 2.2]{SGSamko_AAKilbas_OIMarichev_1993} (see also \cite[Lemma~2.12]{KDiethelm_2010}). Property ($B.5$) is a consequence of ($B.4$) and ($A.1$).
\proofend

\begin{defi}[see {\cite[(2.4.1)]{AAKilbas_HMSrivastava_JJTrujillo_2006}}]
    For a function $x: [0, T] \rightarrow \mathbb{R}^k,$ the (left-sided) Caputo fractional derivative of an order $\alpha \in (0, 1)$ is defined by
    \begin{equation*}
        ({}^C D^\alpha x) (t)
        = \frac{1}{\Gamma(1 - \alpha)} \frac{d}{dt} \int_{0}^{t} \frac{x(\tau) - x(0)}{(t - \tau)^{\alpha}} d\tau, \quad t \in [0, T].
    \end{equation*}
\end{defi}

From the definitions it follows that, for a function $x:[0, T] \rightarrow \mathbb{R}^k,$ if $x(0) = 0,$ then the Caputo and the R.-L. fractional derivatives coincide.

\section{Differential Equation of Fractional Order}

\setcounter{section}{2}
\setcounter{equation}{0}\setcounter{theorem}{0}

Let $n \in \mathbb{N},$ $\alpha \in (0, 1)$ and $T > 0$ be fixed. Let us consider the following Cauchy problem for the ordinary fractional differential equation with the Caputo derivative of the order $\alpha$
\begin{equation} \label{system}
    (^C D^\alpha x) (t) = f(t, x(t)), \quad t \in [0, T], \quad x \in \mathbb{R}^n,
\end{equation}
with the initial condition
\begin{equation}\label{initial_condition}
    x(0) = x_0, \quad x_0 \in \mathbb{R}^n.
\end{equation}
Let the function $f: [0, T] \times \mathbb{R}^n \rightarrow \mathbb{R}^n$ satisfy the following conditions:
\begin{itemize}
  \item[($f.1$)] For any $x \in \mathbb{R}^n,$ the function $f(\cdot, x)$ is measurable on $[0, T].$

  \item[($f.2$)] For any $r \geqslant 0,$ there exists $\lambda_f > 0$ such that
        \begin{equation*}
            \|f(t, x) - f(t, y)\| \leqslant \lambda_f \|x - y\|,
            \quad t \in [0, T], \quad x, y \in B(r).
        \end{equation*}

  \item[($f.3$)] There exists $c_f > 0$ such that
        \begin{equation*}
            \|f(t, x)\| \leqslant (1 + \|x\|) c_f,
            \quad t \in [0, T], \quad x \in \mathbb{R}^n.
        \end{equation*}
\end{itemize}

\begin{defi} \label{Def_solution}
    A function $x: [0, T] \rightarrow \mathbb{R}^n$ is called a solution of Cauchy problem (\ref{system}), (\ref{initial_condition}) if $x(\cdot) \in \{x_0\} + I^\alpha (\Linf([0,T], \mathbb{R}^n))$ and equality (\ref{system}) holds for almost every $t \in [0, T].$
\end{defi}
Here the inclusion $x(\cdot) \in \{x_0\} + I^\alpha (\Linf([0,T], \mathbb{R}^n))$ means that there is a function $y(\cdot) \in I^\alpha (\Linf([0,T], \mathbb{R}^n))$ such that $x(t) = x_0 + y(t),$ $t \in [0, T].$ Note that due to ($A.1$) we have $y(0) = 0,$ and consequently, $x(0) = x_0.$ Therefore, for a function $x(\cdot) \in \{x_0\} + I^\alpha (\Linf([0,T], \mathbb{R}^n)),$ initial condition (\ref{initial_condition}) is automatically satisfied.

\begin{theo} \label{Thm_existence}
    For any initial value $x_0 \in \mathbb{R}^n,$ there exists the unique solution of Cauchy problem (\ref{system}), (\ref{initial_condition}).
\end{theo}
\proof
    By the scheme of the proof from \cite[Lemma~6.2]{KDiethelm_2010} one can show that a function $x: [0, T] \rightarrow \mathbb{R}^n$ is a solution of Cauchy problem (\ref{system}), (\ref{initial_condition}) if and only if $x(\cdot) \in \C([0, T], \mathbb{R}^n)$ and it satisfies the integral equation
    \begin{equation} \label{integral_equation}
        x(t)
        = x_0 + \frac{1}{\Gamma(\alpha)} \int_{0}^{t} \frac{f(\tau, x(\tau))}{(t - \tau)^{1 - \alpha}} d \tau,
        \quad t \in [0, T].
    \end{equation}
    Consequently, it is sufficient to prove the existence and uniqueness of a continuous solution of integral equation (\ref{integral_equation}).

    Let a mapping $F: \C([0, T], \mathbb{R}^n) \rightarrow \C([0, T], \mathbb{R}^n)$ be defined by
    \begin{equation*}
           (F x) (t)
            = x_0 + \frac{1}{\Gamma(\alpha)} \int_{0}^{t} \frac{f(\tau, x(\tau))}{(t - \tau)^{1 - \alpha}} d \tau, 
            \quad t \in [0, T], \quad x(\cdot) \in \C([0, T], \mathbb{R}^n).
    \end{equation*}
    Note that, for any $x(\cdot) \in \C([0, T], \mathbb{R}^n),$ due to ($f.1$)--($f.3$) the function $\varphi(t) = f(t, x(t)),$ $t \in [0, T],$ satisfies the inclusion $\varphi(\cdot) \in \Linf([0, T], \mathbb{R}^n).$ Hence, by ($A.1$) and ($A.2$), the value $(F x) (t)$ is well defined for any $t \in [0, T],$ and $(F x)(\cdot) \in \C([0, T], \mathbb{R}^n).$ Therefore, the definition of $F$ is correct.

    Since a function $x(\cdot) \in \C([0, T], \mathbb{R}^n)$ satisfies integral equation (\ref{integral_equation}) if and only if it is a fixed point of the mapping $F,$ it is sufficient to show the existence and uniqueness of such a fixed point. The proof of this fact is quite standard and follows the scheme described, e.g., in \cite[Theorem~3.1]{JWang_YZhou_2011}. Firstly, due to ($f.2$) one can show that $F$ is continuous. Secondly, the compactness of $F$ follows from ($f.3$) and ($A.3$). Finally, by ($f.3$) and Lemma~\ref{Lem_Bellman_Gronwall}, there exists $r > 0$ such that, for any $x(\cdot) \in \C([0, T], \mathbb{R}^n)$ satisfying $x(t) = \gamma (F x)(t),$ $t \in [0, T],$ with some $\gamma \in (0, 1),$ the inequality $\|x(\cdot)\|_\infty \leqslant r$ is valid. Therefore, by Leray-Shauder theorem (see, e.g., \cite[Theorem~6.2]{EZeidler_1986}), the mapping $F$ has a fixed point.
    Its uniqueness can be shown by the standard argument basing on ($f.2$) and Lemma~\ref{Lem_Bellman_Gronwall}.
\proofend

Let us give some properties of the solution of Cauchy problem (\ref{system}),~(\ref{initial_condition}).
\begin{prop} \label{Prop_motion_properties}
    For any $R_0 > 0,$ there exist $R > 0$ and $H > 0$ such that, for any initial value $x_0 \in B(R_0),$ the solution $x(\cdot)$ of Cauchy problem (\ref{system}), (\ref{initial_condition}) satisfies the inequalities below:
    \begin{equation*}
        \|x(t)\| \leqslant R, \quad
        \|x(t) - x(\tau)\| \leqslant H |t - \tau|^\alpha, \quad t, \tau \in [0, T].
    \end{equation*}
\end{prop}

\proof
    Let $R_0 > 0$ and $H_\infty$ be the constant from ($A.2$). Let us define $R = (1 + R_0) E_\alpha (c_f T^\alpha) - 1,$ $H = H_\infty (1 + R) c_f.$ Let $x_0 \in B(R_0)$ and $x(\cdot)$ be the solution of (\ref{system}), (\ref{initial_condition}). By (\ref{integral_equation}) and ($f.3$), for any $t \in [0, T],$ we have
    \begin{equation*}
        \|x(t)\|
        \leqslant \|x_0\| + \frac{1}{\Gamma(\alpha)} \int_{0}^{t} \frac{\|f(\tau, x(\tau))\|}{(t - \tau)^{1 - \alpha}} d \tau
        \leqslant R_0 + \frac{c_f}{\Gamma(\alpha)} \int_{0}^{t} \frac{1 + \|x(\tau))\|}{(t - \tau)^{1 - \alpha}} d \tau,
    \end{equation*}
    and therefore, according to Lemma~\ref{Lem_Bellman_Gronwall} we obtain
    \begin{equation*}
        \|x(t)\| \leqslant (1 + R_0) E_\alpha (c_f T^\alpha) - 1 = R, \quad t \in [0, T].
    \end{equation*}
    Further, from ($f.3$) it follows that the function $\varphi(t) = f(t, x(t)),$ $t \in [0, T],$ satisfies the inequalities
    \begin{equation*}
        \|\varphi(t)\| \leqslant (1 + \|x(t)\|) c_f \leqslant (1 + R) c_f \text{ for a.e. } t \in [0, T],
    \end{equation*}
    wherefrom, due to (\ref{integral_equation}) and ($A.2$), for any $t, \tau \in [0, T],$ we derive
    \begin{equation*}
        \|x(t) - x(\tau)\|
        = \|(I^\alpha \varphi)(t) - (I^\alpha \varphi)(\tau)\|
        \leqslant H_\infty (1 + R) c_f |t - \tau|^\alpha
        = H |t - \tau|^\alpha.
    \end{equation*}
    The proposition is proved.
\proofend

\section{Fractional Derivative of a Convex Lyapunov Function}

\setcounter{section}{3}
\setcounter{equation}{0}\setcounter{theorem}{0}

Let a function $V: \mathbb{R}^n \rightarrow \mathbb{R}$ satisfy the following conditions:
\begin{itemize}
  \item[($V.1$)] The function $V(\cdot)$ is convex on $\mathbb{R}^n$ and $V(0) = 0.$

  \item[($V.2$)] The function $V(\cdot)$ is differentiable (and therefore continuous) on~$\mathbb{R}^n.$

  \item[($V.3$)] For any $r \geqslant 0,$ there exists $\lambda_V > 0$ such that
    \begin{equation*}
        \|\nabla V(x) - \nabla V(y)\| \leqslant \lambda_V \|x - y\|, \quad x, y \in B(r),
    \end{equation*}
    where $\nabla V(\cdot)$ is the gradient of the function $V(\cdot).$
\end{itemize}

According to \cite[Theorem~1]{WChen_HDai_YSong_ZZhang_2017}, for a sufficiently smooth function $x: [0, T] \rightarrow \mathbb{R}^n,$ $x(0) = 0,$ if we denote $y(t) = V(x(t)),$ $t \in [0, T],$ then, for any $t \in [0, T],$ the following inequality holds:
\begin{equation} \label{Prop_convex_smooth_main_1}
    (D^\alpha y) (t) \leqslant \langle \nabla V(x(t)), (D^\alpha x) (t) \rangle.
\end{equation}
The proof of this fact is based on representation formula (\ref{representation_formula}) (see Proposition~\ref{Prop_convex_smooth} below). Therefore, in particular, it substantially uses differentiability properties of the function $x(\cdot).$ However, the solution of Cauchy problem (\ref{system}), (\ref{initial_condition}) may be nowhere differentiable (see, e.g., \cite{BRoss_SGSamko_ERLove_1994}). Hence, the technique used in the proof can not be directly applied to prove inequality (\ref{Prop_convex_smooth_main_1}) for the case when $x(\cdot)$ is the solution of Cauchy problem (\ref{system}), (\ref{initial_condition}).

The goal of this section is to establish estimate (\ref{Prop_convex_smooth_main_1}) for any function $x(\cdot) \in I^\alpha(\Linf([0, T], \mathbb{R}^n)).$ The proof is carried out in several stages. Firstly, the smooth case, when $x(\cdot) \in \Lip^0([0, T], \mathbb{R}^n),$ is studied. After that, it is proved that any function from $I^\alpha(\Linf([0, T], \mathbb{R}^n))$ can be approximated by functions from $\Lip^0([0, T], \mathbb{R}^n)$ with the uniformly bounded derivatives of the order $\alpha.$ Finally, in the general case, applying for the smooth approximating functions results that have been already obtained, the estimate (\ref{Prop_convex_smooth_main_1}) is proved for any function $x(\cdot) \in I^\alpha(\Linf([0, T], \mathbb{R}^n)).$

\subsection{Smooth Case}

\begin{prop} \label{Prop_convex_smooth}
    Let $x(\cdot) \in \Lip^0([0, T], \mathbb{R}^n)$ and $y(t) = V(x(t)),$ $t \in [0, T].$ Then the inclusion $y(\cdot) \in \Lip^0 ([0,T], \mathbb{R})$ is valid and inequality (\ref{Prop_convex_smooth_main_1}) holds for every $t \in [0, T].$ Moreover, for any $r \geqslant 0$ and $w \geqslant 0,$ there exists $a \geqslant 0$ such that, for any $x(\cdot) \in \Lip^0([0, T], \mathbb{R}^n),$ if
    \begin{equation}\label{r_s}
        \|x(\cdot)\|_\infty \leqslant r, \quad \|(D^\alpha x)(\cdot)\|_\infty \leqslant w,
    \end{equation}
    then the function $y(t) = V(x(t)),$ $t \in [0, T],$ satisfies the inequality
    \begin{equation}\label{Prop_convex_smooth_main_2}
        \|(D^\alpha y)(\cdot)\|_\infty \leqslant a.
    \end{equation}
\end{prop}
\proof
    Let $x(\cdot) \in \Lip^0([0, T], \mathbb{R}^n)$ and $L \geqslant 0$ be the Lipschitz constant of $x(\cdot).$ Let $r \geqslant 0$ and $w \geqslant 0$ satisfy inequalities (\ref{r_s}). Due to ($V.3$), by the number $r,$ let us choose $\lambda_V > 0$ and put $M_V = \max_{x \in B(r)} \|\nabla V(x)\|.$ Let $H_\infty$ be the constant from ($A.2$). Then from ($B.2$) it follows that the function $x(\cdot)$ is H\"{o}lder continuous of the order $\alpha$ with the constant $H = H_\infty w.$

    Let $y(t) = V(x(t)),$ $t \in [0, T].$ Let us show that $y(\cdot) \in \Lip^0([0, T], \mathbb{R}^n).$ From ($V.1$) it follows that $y(0) = 0.$ Further, let $t, \tau \in [0, T].$ Due to ($V.2$), by the mean value theorem, there exists $\gamma \in [0, 1]$ such that, for the vector $z = \gamma x(t) + (1 - \gamma) x(\tau),$ we have
    \begin{equation} \label{y_t_y_tau}
        y(t) - y(\tau) = V(x(t)) - V(x(\tau)) = \langle \nabla V(z), x(t) - x(\tau) \rangle.
    \end{equation}
    Hence, since $z \in B(r),$ by the choice of $L$ and $M_V,$ we obtain
    \begin{equation*}
        |y(t) - y(\tau)| \leqslant \|\nabla V(z)\| \|x(t) - x(\tau)\| \leqslant M_V L |t - \tau|.
    \end{equation*}
    Therefore, the function $y(\cdot)$ is Lipschitz continuous.

    The proof of inequality (\ref{Prop_convex_smooth_main_1}) follows the scheme from \cite[Theorem~1]{WChen_HDai_YSong_ZZhang_2017}. But it seems convenient to give this proof because its main part is used in the proof of the last part of the proposition.

    Since $x(\cdot), y(\cdot) \in \Lip^0([0, T], \mathbb{R}^n),$ then inequality (\ref{Prop_convex_smooth_main_1}) for $t = 0$ follows from ($B.5$). Let $t \in (0, T].$ Due to ($V.2$), by the chain rule, we have $\dot{y}(t) = \langle \nabla V(x(t)), \dot{x}(t) \rangle$ for almost every $t \in [0, T].$ Therefore, by ($B.4$), inequality (\ref{Prop_convex_smooth_main_1}) multiplied by $\Gamma(1 - \alpha)$ can be rewritten as follows
    \begin{equation} \label{Prop_convex_smooth_1}
        \int_{0}^{t} \frac{\langle \nabla V(x(\tau)), \dot{x}(\tau) \rangle}{(t - \tau)^\alpha} d \tau
        \leqslant  \int_{0}^{t} \frac{\langle \nabla V(x(t)), \dot{x}(\tau) \rangle}{(t - \tau)^\alpha} d \tau.
    \end{equation}
    Let us consider the function
    \begin{equation*}
        \varphi(\tau) = V(x(\tau)) - V(x(t)) - \langle \nabla V(x(t)), x(\tau) - x(t) \rangle, \quad \tau \in [0, t].
    \end{equation*}
    Then $\varphi(\cdot) \in \Lip([0, t], \mathbb{R})$ and
    \begin{equation*}
        \dot{\varphi}(\tau)
        = \langle \nabla V(x(\tau)) - \nabla V(x(t)), \dot{x}(\tau) \rangle
        \text{ for a.e. } \tau \in [0, t].
    \end{equation*}
    Hence,
    \begin{equation*}
        \int_{0}^{t} \frac{\langle \nabla V(x(\tau)) - \nabla V(x(t)), \dot{x}(\tau) \rangle}{(t - \tau)^\alpha} d \tau
        = \int_{0}^{t} \frac{\dot{\varphi}(\tau)}{(t - \tau)^\alpha} d \tau,
    \end{equation*}
    and in order to prove inequality (\ref{Prop_convex_smooth_1}) it is sufficient to show that
    \begin{equation} \label{Prop_convex_smooth_2}
        \int_{0}^{t} \frac{\dot{\varphi}(\tau)}{(t - \tau)^\alpha} d \tau \leqslant 0.
    \end{equation}

    Let us prove that
    \begin{equation} \label{Prop_convex_smooth_3}
        0 \leqslant \varphi(\tau) \leqslant \lambda_V H^2 (t - \tau)^{2 \alpha}, \quad \tau \in [0, t].
    \end{equation}
    Let $\tau \in [0, t].$ Let $\gamma \in [0, 1]$ and $z = \gamma x(t) - (1 - \gamma) x(\tau) \in B(r)$ be such that (\ref{y_t_y_tau}) is valid. Therefore, we have
    \begin{equation*}
       \varphi(\tau)
        = \langle \nabla V(z), x(\tau) - x(t) \rangle - \langle \nabla V(x(t)), x(\tau) - x(t) \rangle.
    \end{equation*}
    Consequently, by the choice of $\lambda_V$ and $H,$ we obtain
    \begin{equation*}
        \begin{array}{c}
            \varphi(\tau)
            \leqslant \| \nabla V(z) - \nabla V(x(t))\| \|x(\tau) - x(t) \|
            \leqslant \lambda_V \|z - x(t)\| \|x(\tau) - x(t)\| \\[0.5em]
            \leqslant \lambda_V \|x(\tau) - x(t)\|^2
            \leqslant \lambda_V H (t - \tau)^{2 \alpha}.
        \end{array}
    \end{equation*}
    On the other hand, due to ($V.1$) and ($V.2$), by the differentiation of convex functions theorem (see, e.g., \cite[Theorem~25.1]{RTRockafellar_1972}), we have
    \begin{equation*}
        V(x(\tau)) - V(x(t)) \geqslant \langle \nabla V(x(t)), x(\tau) - x(t) \rangle,
    \end{equation*}
    and hence,
    \begin{equation*}
        \varphi(\tau)
        = V(x(\tau)) - V(x(t)) - \langle \nabla V(x(t)), x(\tau) - x(t) \rangle
        \geqslant 0.
    \end{equation*}

    Taking (\ref{Prop_convex_smooth_3}) into account, by the integration by parts formula, we derive
    \begin{equation} \label{Prop_convex_smooth_4}
        \int_{0}^t \frac{\dot{\varphi}(\tau)}{(t - \tau)^\alpha} d \tau
        = - \frac{\varphi(0)}{t^\alpha} - \alpha \int_{0}^t \frac{\varphi(\tau)}{(t - \tau)^{\alpha + 1}} d \tau.
    \end{equation}
    Thus, inequality (\ref{Prop_convex_smooth_2}) follows from (\ref{Prop_convex_smooth_3}) and (\ref{Prop_convex_smooth_4}).

    Let us prove the remaining part of the proposition. Let $r \geqslant 0$ and $w \geqslant 0.$ Let us define
    \begin{equation} \label{A}
        a = 2 \lambda_V H^2 T^\alpha/\Gamma(1 - \alpha) + M_V w.
    \end{equation}
    Let $x(\cdot) \in \Lip^0([0, T], \mathbb{R}^n)$ satisfy inequalities (\ref{r_s}) and $y(t) = V(x(t)),$ $t \in [0, T].$ Let us show that inequality (\ref{Prop_convex_smooth_main_2}) is valid with this number $a.$

    If $t = 0,$ then inequality (\ref{Prop_convex_smooth_main_2}) follows from ($B.5$). Let $t \in (0, T].$ By analogy with the previous arguments, we have
    \begin{equation} \label{A_1}
        (D^\alpha y) (t)
        = \frac{1}{\Gamma(1 - \alpha)} \int_{0}^{t} \frac{\dot{\varphi}(\tau)}{(t - \tau)^\alpha} d \tau
        + \langle \nabla V(x(t)), (D^\alpha x) (t) \rangle.
    \end{equation}
    From (\ref{Prop_convex_smooth_3}) and (\ref{Prop_convex_smooth_4}) we derive
    \begin{equation} \label{A_2}
        \begin{array}{c}
            \displaystyle
            0 \geqslant \int_{0}^{t} \frac{\dot{\varphi}(\tau)}{(t - \tau)^\alpha} d \tau
            \geqslant - \frac{\lambda_V H^2 t^{2 \alpha}}{t^\alpha}
            - \alpha \int_{0}^t \frac{\lambda_V H^2 (t - \tau)^{2 \alpha}}{(t - \tau)^{\alpha + 1}} d \tau \\[1em]
            = - 2 \lambda_V H^2 t^\alpha
            \geqslant - 2 \lambda_V H^2 T^\alpha,
        \end{array}
    \end{equation}
    and due to the choice of $M_V$ we obtain
    \begin{equation} \label{A_3}
        |\langle \nabla V(x(t)), (D^\alpha x) (t) \rangle|
        \leqslant \|\nabla V(x(t))\| \|(D^\alpha x) (t)\|
        \leqslant M_V w.
    \end{equation}
    Thus, inequality (\ref{Prop_convex_smooth_main_2}) with $a$ defined in (\ref{A}) follows from (\ref{A_1})--(\ref{A_3}).
\proofend

\subsection{Approximation}

\begin{prop} \label{Prop_approximation_L_infty}
    Let $\varphi(\cdot) \in \Linf([0, T], \mathbb{R}^n)$ and $p \in [1, \infty).$ Then, for any $\varepsilon > 0,$ there exists $\overline{\varphi}(\cdot) \in \Lip^0 ([0, T], \mathbb{R}^n)$ such that $\|\overline{\varphi}(\cdot)\|_\infty \leqslant \sqrt{n} \|\varphi(\cdot)\|_\infty$ and $\|\varphi(\cdot) - \overline{\varphi}(\cdot)\|_p \leqslant \varepsilon.$
\end{prop}
\proof
    Let $\varphi(\cdot) \in \Linf([0, T], \mathbb{R}^n),$ $p \in [1,\infty)$ and $\varepsilon > 0.$ Let $\xi > 0$ be such that $(1 + \sqrt{n}) \|\varphi(\cdot)\|_\infty \ \xi^{1/p} \leqslant \varepsilon/2.$ Applying Lusin theorem (see, e.g., \cite[Theorem~2.24]{WRudin_1987}) to each coordinate of $\varphi(\cdot),$ one can find a function $\psi(\cdot) \in \C([0, T], \mathbb{R}^n)$ such that the set $E = \{t \in [0, T]: \varphi (t) \neq \psi (t)\}$ has measure less than $\xi$ and $\|\psi(\cdot)\|_\infty \leqslant \sqrt{n} \|\varphi(\cdot)\|_\infty.$ Since
    \begin{equation*}
        \begin{array}{c}
            \displaystyle
            \|\psi(\cdot) - \varphi(\cdot)\|_p^p
            = \int_{0}^{T} \|\psi (t) - \varphi (t)\|^p dt
            = \int_{E} \|\psi (t) - \varphi (t)\|^p dt \\[1em]
            \leqslant \|\psi (\cdot) - \varphi (\cdot)\|_\infty^p \ \xi
            \leqslant (\sqrt{n} + 1)^p \|\varphi(\cdot)\|_\infty^p \ \xi,
        \end{array}
    \end{equation*}
    then, by the choice of $\xi$ we have
    \begin{equation} \label{psi-varphi}
        \|\psi(\cdot) - \varphi(\cdot)\|_p \leqslant (\sqrt{n} + 1) \|\varphi(\cdot)\|_\infty \ \xi^{1/p} \leqslant \varepsilon/2.
    \end{equation}

    Further, let $\eta > 0$ and $(1 + T)^{1/p} \eta \leqslant \varepsilon/2.$ Since $\psi(\cdot) \in \C([0, T], \mathbb{R}^n),$ one can choose $\delta_1 > 0$ such that, for any $t, \tau \in [0, T],$ if $|t - \tau| \leqslant \delta_1,$ then $\|\psi(t) - \psi(\tau)\| \leqslant \eta/2.$ Let $\delta_2 > 0$ and $2 \sqrt{n} \|\varphi(\cdot)\|_\infty \delta_2^{1/p} \leqslant \eta.$ Let $\delta = \min\{\delta_1, \delta_2\}$ and $N \in \mathbb{N}$ satisfy the inequality $T/N \leqslant \delta.$ Let us denote $t_i = T i/ N,$ $i \in \overline{0, N},$ and define a piecewise linear function $\overline{\varphi}(\cdot) \in \Lip^0([0, T], \mathbb{R}^n):$
    \begin{equation*}
        \begin{array}{c}
            \displaystyle
            \overline{\varphi} (t) = \psi(t_1) t / \delta,
            \quad t \in [t_0, t_1], \\[0.5em]
            \displaystyle
            \overline{\varphi} (t) = \psi(t_i)
            + (\psi(t_{i+1}) - \psi(t_i)) (t - t_i) / \delta,
            \ \ t \in [t_i, t_{i+1}], \ \ i \in \overline{1, N - 1}.
        \end{array}
    \end{equation*}
    From the definition it follows that
    \begin{equation*}
        \|\overline{\varphi}(\cdot)\|_\infty
        = \max_{i \in \overline{1, N}} \|\psi(t_i)\|
        \leqslant \|\psi(\cdot)\|_\infty \leqslant \sqrt{n} \|\varphi(\cdot)\|_\infty.
    \end{equation*}
    For $t \in [t_0, t_1],$ we obtain
    \begin{equation*}
        \|\overline{\varphi}(t) - \psi (t)\|
        = \| \psi(t_1) t / \delta - \psi (t)\|
        \leqslant 2 \|\psi(\cdot)\|_\infty \leqslant 2 \sqrt{n} \|\varphi(\cdot)\|_\infty,
    \end{equation*}
    and for $t \in [t_i, t_{i+1}],$ $i \in \overline{1, N-1},$ according to the choice of $\delta_1$ we derive
    \begin{equation*}
        \begin{array}{c}
            \|\overline{\varphi} (t) - \psi (t)\|
            \leqslant \| \psi(t_i) - \psi (t) \|
            + \|(\psi(t_{i + 1}) - \psi(t_i) (t - t_i) / \delta \| \\[0.5em]
            \leqslant \| \psi(t_i) - \psi (t) \|
            + \|\psi(t_{i + 1}) - \psi(t_i)\|
            \leqslant \eta.
        \end{array}
    \end{equation*}
    Consequently, due to the choice of $\delta_2$ and $\eta$ we have
    \begin{equation*}
        \begin{array}{c}
            \displaystyle
            \|\overline{\varphi} (\cdot) - \psi(\cdot) \|_p^p
            = \int_{t_0}^{t_1} \|\overline{\varphi} (t) - \psi (t) \|^p dt
            + \int_{t_1}^{t_N} \|\overline{\varphi} (t) - \psi (t) \|^p dt \\[1em]
            \displaystyle
            \leqslant (2 \sqrt{n} \|\varphi(\cdot)\|_\infty)^p \delta + T \eta^p
            \leqslant (1 + T) \eta^p
            \leqslant \varepsilon^p/2^p.
        \end{array}
    \end{equation*}
    Therefore,
    \begin{equation} \label{varphi-varphi}
        \|\overline{\varphi}(\cdot) - \psi(\cdot) \|_p \leqslant \varepsilon/2.
    \end{equation}

    Thus, from (\ref{psi-varphi}) and (\ref{varphi-varphi}) it follows that $\|\varphi(\cdot) - \overline{\varphi}(\cdot)\|_p \leqslant \varepsilon.$
\proofend

\begin{cor} \label{Cor_approximation_I_a+^alpha}
    For any $x(\cdot) \in I^\alpha (\Linf([0, T], \mathbb{R}^n))$ and $p \in (1/\alpha, \infty),$ there exists a sequence $\{x_k(\cdot)\}_{k=1}^\infty \subset \Lip^0 ([0, T], \mathbb{R}^n)$ such that the inequality $\|(D^\alpha x_k)(\cdot)\|_\infty \leqslant \sqrt{n} \|(D^\alpha x)(\cdot)\|_\infty$ is valid for any $k \in \mathbb{N}$ and
    \begin{equation} \label{Cor_approximation_I_a+^alpha_main_2}
        \lim_{k \rightarrow \infty} \|x_k(\cdot) - x(\cdot)\|_\infty = 0, \quad
        \lim_{k \rightarrow \infty} \|(D^\alpha x_k)(\cdot) - (D^\alpha x)(\cdot)\|_p = 0.
    \end{equation}
\end{cor}
\proof
     Let $x(\cdot) \in I^\alpha (\Linf([0, T], \mathbb{R}^n))$ and $p \in (1/\alpha, \infty).$ For the function $\varphi(t) = (D^\alpha x)(t),$ $t \in [0, T],$ by Proposition~\ref{Prop_approximation_L_infty}, for every $k \in \mathbb{N},$ one can choose $\varphi_k(\cdot) \in \Lip^0([0, T], \mathbb{R}^n)$ such that $\|\varphi_k(\cdot)\|_\infty \leqslant \sqrt{n} \|\varphi(\cdot)\|_\infty$ and $\|\varphi(\cdot) - \varphi_k(\cdot)\|_p \leqslant 1/k.$ Therefore, $\|\varphi(\cdot) - \varphi_k(\cdot)\|_p \rightarrow 0$ when $k \rightarrow \infty.$ Let $x_k(t) = (I^\alpha \varphi_k)(t),$ $t \in [0, T],$ $k \in \mathbb{N}.$ Due to ($A.4$) we obtain $x_k (\cdot) \in \Lip^0([0, T], \mathbb{R}^n),$ $k \in \mathbb{N}.$ By ($A.3$), we have $\|x_k(\cdot) - x(\cdot)\|_\infty \rightarrow 0$ when $k \rightarrow \infty,$ and consequently, the first relation in (\ref{Cor_approximation_I_a+^alpha_main_2}) is valid. For every $k \in \mathbb{N},$ according to ($B.3$) we get $(D^\alpha x_k)(t) = \varphi_k(t)$ for almost every $t \in [0, T],$ and therefore, the second relation in (\ref{Cor_approximation_I_a+^alpha_main_2}) holds. The corollary is proved.
\proofend

\subsection{General Case}

\begin{lem} \label{Lem_differentiating_of_convex}
    Let $x(\cdot) \in I^\alpha (\Linf([0, T], \mathbb{R}^n))$ and $y(t) = V(x(t)),$ $t \in [0, T].$ Then the inclusion $y(\cdot) \in I^\alpha(\Linf([0, T], \mathbb{R}))$ is valid and inequality (\ref{Prop_convex_smooth_main_1}) holds for almost every $t \in [0, T].$
\end{lem}
\proof
    Let $x(\cdot) \in I^\alpha (\Linf([0, T], \mathbb{R}^n))$ and $p \in (1/\alpha, \infty).$ Due to Corollary~\ref{Cor_approximation_I_a+^alpha} one can choose a sequence $\{x_k(\cdot)\}_{k = 1}^\infty \subset \Lip^0([0, T], \mathbb{R}^n)$ such that $\|\varphi_k(\cdot)\|_\infty \leqslant \sqrt{n} \|\varphi(\cdot)\|_\infty,$ $k \in \mathbb{N},$ and
    \begin{equation} \label{Lem_differentiating_of_convex_convergence_x_varphi}
        \lim_{k \rightarrow \infty} \|x_k(\cdot) - x(\cdot)\|_\infty = 0, \quad
        \lim_{k \rightarrow \infty} \|\varphi_k(\cdot) - \varphi(\cdot)\|_p = 0,
    \end{equation}
    where $\varphi_k(t) = (D^\alpha x_k)(t),$ $t \in [0, T],$ $k \in \mathbb{N},$ and $\varphi(t) = (D^\alpha x)(t),$ $t \in [0, T].$ Note that, by ($B.1$) and ($B.5$), we have $\varphi(\cdot), \varphi_k(\cdot) \in \Linf ([0, T], \mathbb{R}^n),$ $k \in \mathbb{N}.$

    Let $r = \sup_{k \in \mathbb{N}} \|x_k(\cdot)\|_\infty$ and $w = \sup_{k \in \mathbb{N}} \|\varphi_k(\cdot)\|_\infty.$ In particular, from the first relation in (\ref{Lem_differentiating_of_convex_convergence_x_varphi}) it follows that $\|x(\cdot)\|_\infty \leqslant r.$ Due to ($V.3$), by the number $r,$ let us choose $\lambda_V > 0$ and put $M_V = \max_{x \in B(r)} \|\nabla V(x)\|.$

    Let $y_k(t) = V(x_k(t)),$ $t \in [0, T],$ $k \in \mathbb{N}.$ For every $k \in \mathbb{N},$ according to Proposition~\ref{Prop_convex_smooth} the inclusion $y_k (\cdot) \in \Lip^0([0, T], \mathbb{R})$ is valid and, for the function $\psi_k(t) = (D^\alpha y_k)(t),$ $t \in [0, T],$ the following inequality holds:
    \begin{equation} \label{Lem_differentiating_of_convex_estimate}
        \psi_k (t) \leqslant \langle \nabla V(x_k(t)), \varphi_k(t) \rangle, \quad t \in [0, T].
    \end{equation}
    Moreover, there exists $a \geqslant 0$ such that $\|\psi_k(\cdot)\|_\infty \leqslant a,$ $k \in \mathbb{N}.$

    Let us consider the set
    \begin{equation*}
        K = \big\{ \psi(\cdot) \in \Lp([0, T], \mathbb{R}): \|\psi(\cdot)\|_\infty \leqslant a\big\}.
    \end{equation*}
    This set is weakly sequentially compact in $\Lp([0, T], \mathbb{R}).$ Indeed, $K$ is convex, bounded and, applying \cite[Theorem~3.12]{WRudin_1987}, one can show that $K$ is closed. Consequently, from \cite[Ch~III, \S~3, Theorem~2]{LVKantorovich_GPAkilov_1982} it follows that $K$ is weakly closed. Therefore, by \cite[Ch.~V, \S~7, Theorem~7]{LVKantorovich_GPAkilov_1982} set $K$ is weakly compact as a weakly closed subset of a weakly compact set. Hence, due to \cite[Ch.~VIII, \S~2, Corollary]{LVKantorovich_GPAkilov_1982} this set is weakly sequentially compact.

    Since $\{\psi_k(\cdot)\}_{k = 1}^\infty \subset K,$ we can assume that the sequence $\{\psi_k(\cdot)\}_{k = 1}^\infty$ converges weakly to a function $\overline{\psi}(\cdot) \in K.$ Note that, $\overline{\psi}(\cdot) \in \Linf([0, T], \mathbb{R}).$ From ($A.3$) and \cite[Proposition~3.3]{JBConway_1985} we obtain $\|(I^\alpha \psi_k)(\cdot) - (I^\alpha \overline{\psi})(\cdot)\|_\infty \rightarrow 0$ when $k \rightarrow \infty.$ Due to ($B.2$) we have $y_k(t) = (I^\alpha \psi_k)(t),$ $t \in [0, T],$ $k \in \mathbb{N},$ therefore, $\|y_k(\cdot) - (I^\alpha \overline{\psi}) (\cdot)\|_\infty \rightarrow 0$ when $k \rightarrow \infty.$ On the other hand, from ($V.2$) and the first relation in (\ref{Lem_differentiating_of_convex_convergence_x_varphi}) it follows that $\|y_k(\cdot) - y(\cdot)\|_\infty \rightarrow 0$ when $k \rightarrow \infty.$ Consequently, $y(t) = (I^\alpha \overline{\psi})(t),$ $t \in [0, T].$ Hence, $y(\cdot) \in I^\alpha (\Linf([0, T], \mathbb{R}))$ and due to ($B.3$) we have $\overline{\psi}(t) = (D^\alpha y)(t)$ for almost every $t \in [0, T].$

    Let us prove that inequality (\ref{Prop_convex_smooth_main_1}) holds for almost every $t \in [0, T].$ Let $j \in \mathbb{N}.$ Since the sequence $\{\psi_k\}_{k = j}^\infty$ converges weakly to $\overline{\psi}(\cdot),$ then, by \cite[Theorem~3.13]{WRudin_1991}, there exists a convex combination $\xi_j(t) = \sum_{i = 1}^{n_j} \alpha_{ij} \psi_{k_{ij}}(t),$ $t \in [0, T],$ that satisfies the inequality $\|\xi_j(\cdot) - \overline{\psi}(\cdot)\|_p \leqslant 1/j.$ Here $n_j \in \mathbb{N},$ $k_{ij} \in \mathbb{N},$ $k_{ij} \geqslant j,$ $\alpha_{ij} \in [0, 1],$ $i \in \overline{1, n_j},$ and $\sum_{i = 1}^{n_j} \alpha_{ij} = 1.$ Thus, for the sequence $\{\xi_j(\cdot)\}_{j = 1}^\infty \in \Linf([0, T], \mathbb{R})$ we have $\|\xi_j(\cdot) - \overline{\psi}(\cdot)\|_p \rightarrow 0,$ $j \rightarrow \infty.$

    Let us consider the following functions:
    \begin{equation*}
        \begin{array}{c}
            z(t) = \langle \nabla V(x(t)), \varphi(t) \rangle, \quad t \in [0, T], \\
            \displaystyle
            z_j(t) = \sum_{i = 1}^{n_j} \alpha_{ij} \langle \nabla V(x_{k_{ij}}(t)), \varphi_{k_{ij}}(t) \rangle,
            \quad t \in [0, T], \quad j \in \mathbb{N}.
        \end{array}
    \end{equation*}
    Note that $z(\cdot), z_j(\cdot) \in \Linf([0, T], \mathbb{R}),$ $j \in \mathbb{N}.$ Let us show that $\|z_j(\cdot) - z(\cdot)\|_p \rightarrow 0$ when $j \rightarrow \infty.$ Let $\varepsilon > 0.$ In accordance with (\ref{Lem_differentiating_of_convex_convergence_x_varphi}) there exists $J > 0$ such that, for any $j \in \mathbb{N},$ $j \geqslant J,$ the inequality below is valid:
    \begin{equation*}
        w \lambda_V T^{1/p} \sup_{k \geqslant j} \|x_k(\cdot) - x(\cdot)\|_\infty
        + M_V \sup_{k \geqslant j} \|\varphi_k(\cdot) - \varphi(\cdot) \|_p
        \leqslant \varepsilon.
    \end{equation*}
    Let $j \in \mathbb{N}$ and $j \geqslant J.$ For almost every $t \in [0, T],$ according to the choice of $w,$ $\lambda_V$ and $M_V$ we derive
    \begin{equation*}
        \begin{array}{c}
            \displaystyle
            |z_j(t) - z(t)|
            \leqslant \sum_{i = 1}^{n_j} \alpha_{ij} \|\nabla V(x_{k_{ij}}(t)) - \nabla V(x(t))\| \|\varphi_{k_{ij}}(t) \| \\
            \displaystyle
            + \sum_{i = 1}^{n_j} \alpha_{ij} \|\nabla V(x(t))\| \|\varphi_{k_{ij}}(t) - \varphi(t) \| \\
            \displaystyle
            \leqslant w \lambda_V \sum_{i = 1}^{n_j} \alpha_{ij} \|x_{k_{ij}}(\cdot) - x(\cdot)\|_\infty
            + M_V \sum_{i = 1}^{n_j} \alpha_{ij} \|\varphi_{k_{ij}}(t) - \varphi(t) \|.
        \end{array}
    \end{equation*}
    We have $k_{ij} \geqslant j \geqslant J,$ $i \in \overline{1, n_j},$ and hence, due to the choice of $J$ we obtain
    \begin{equation*}
        \|z_j(\cdot) - z(\cdot)\|_p
        \leqslant \big\| w \lambda_V \sup_{k \geqslant j} \|x_k(\cdot) - x(\cdot)\|_\infty \big\|_p
        + M_V \sup_{k \geqslant j} \|\varphi_{k}(\cdot) - \varphi(\cdot)\|_p \\[0.5em]
        \leqslant \varepsilon.
    \end{equation*}

    For any $t \in [0, T]$ and $j \in \mathbb{N},$ from (\ref{Lem_differentiating_of_convex_estimate}) it follows that
    \begin{equation} \label{xi_leq_z}
        \xi_j (t)
        = \sum_{i = 1}^{n_j} \alpha_{ij} \psi_{k_{ij}} (t)
        \leqslant \sum_{i = 1}^{n_j} \alpha_{ij} \langle \nabla V(x_{k_{ij}}(t)), \varphi_{k_{ij}}(t) \rangle
        = z_j(t).
    \end{equation}
    Since $\|\xi_j(\cdot) - \overline{\psi}(\cdot)\|_p \rightarrow 0$ and $\|z_j(\cdot) - z(\cdot)\|_p \rightarrow 0$ when $j \rightarrow \infty,$ then due to \cite[Theorem~3.12]{WRudin_1987} we can assume that $|\xi_j(t) - \overline{\psi}(t)| \rightarrow 0$ and $|z_j(t) - z(t)| \rightarrow 0$ for almost every $t \in [0, T].$ Therefore, for almost every $t \in [0, T],$ letting $j$ to $\infty$ in (\ref{xi_leq_z}), we derive $\overline{\psi}(t) \leqslant z(t).$ Consequently, taking into account that $\overline{\psi}(t) = (D^\alpha y)(t)$ for almost every $t \in [0, T]$ and $z(t) = \langle \nabla V(x(t)), (D^\alpha x)(t) \rangle,$ $t \in [0, T],$ we obtain the validity of inequality (\ref{Prop_convex_smooth_main_1}) for almost every $t \in [0, T].$ The lemma is proved.
\proofend

For the case when $V(x) = \|x\|^2,$ $x \in \mathbb{R}^n,$ we obtain the following result.
\begin{cor} \label{Cor_differentiating_of_square}
    Let $x(\cdot) \in I^\alpha (\Linf([0, T], \mathbb{R}^n))$ and $y(t) = \|x(t)\|^2,$ $t \in [0, T].$ Then the inclusion $y(\cdot) \in I^\alpha(\Linf([0, T], \mathbb{R}))$ is valid and the inequality $(D^\alpha y) (t) \leqslant 2 \langle x(t), (D^\alpha x)(t) \rangle$ holds for almost every $t \in [0, T].$
\end{cor}

\section{Conflict-Controlled Dynamical System of Fractional Order}

\setcounter{section}{4}
\setcounter{equation}{0}\setcounter{theorem}{0}

Let us consider a conflict-controlled dynamical system which motion is described by the fractional differential equation
\begin{equation} \label{system_u_v}
    \begin{array}{c}
        (^C D^\alpha x) (t) = g(t, x(t), u(t), v(t)), \quad t \in [0, T], \\[0.5em]
        x(t) \in \mathbb{R}^n, \quad u(t) \in P \subset \mathbb{R}^{n_u}, \quad v(t) \in Q \subset \mathbb{R}^{n_v},
    \end{array}
\end{equation}
with the initial condition
\begin{equation}\label{initial_condition_u_v}
    x(0) = x_0, \quad x_0 \in \mathbb{R}^n.
\end{equation}

Here $t$ is the time variable, $x$ is the state vector, $u$ is the control vector and $v$ is the vector of unknown disturbances; $n_u, n_v \in \mathbb{N};$ $P$ and $Q$ are compact sets; $x_0$ is the initial value of the state vector. The function $g: [0, T] \times \mathbb{R}^n \times P \times Q \rightarrow \mathbb{R}^n$ satisfies the following conditions:
\begin{itemize}
  \item[($g.1$)] The function $g(\cdot)$ is continuous.
  \item[($g.2$)] For any $r \geqslant 0,$ there exists $\lambda_g > 0$ such that
      \begin{equation*}
        \begin{array}{c}
            \|g(t, x, u, v) - g(t, y, u, v)\| \leqslant \lambda_g \|x - y\|, \\[0.5em]
            t \in [0, T], \quad x, y \in B(r), \quad u \in P, \quad v \in Q.
        \end{array}
      \end{equation*}
  \item[($g.3$)] There exits $c_g > 0$ such that
    \begin{equation*}
        \|g(t, x, u, v)\| \leqslant (1 + \|x\|) c_g,
        \ \ t \in [0, T], \ \ x \in \mathbb{R}^n, \ \ u \in P, \ \ v \in Q.
    \end{equation*}
  \item[($g.4$)] For any $t \in [0, T]$ and $x, s \in \mathbb{R}^n,$ the following equality holds:
    \begin{equation*}
            \displaystyle
            \min_{u \in P} \max_{v \in Q} \langle s, g(t, x, u, v) \rangle
            = \max_{v \in Q} \min_{u \in P} \langle s, g(t, x, u, v) \rangle.
    \end{equation*}
\end{itemize}

It should be noted here that these conditions are quite standard for the differential games theory (see, e.g., \cite[pp.~7, 8]{NNKrasovskii_AISubbotin_1988}).

\begin{defi}
    Admissible control and disturbance realizations are measurable functions $u: [0, T) \rightarrow P$ and $v: [0, T) \rightarrow Q,$ respectively. The corresponding sets of all admissible control $u(\cdot)$ and disturbance $v(\cdot)$ realizations are denoted by $\mathcal{U}$ and $\mathcal{V}.$
\end{defi}

\begin{defi} \label{Def_x}
    A motion of system (\ref{system_u_v}), (\ref{initial_condition_u_v}) that corresponds to an initial value $x_0 \in \mathbb{R}^n$ and realizations $u(\cdot) \in \mathcal{U},$ $v(\cdot) \in \mathcal{V}$ is a solution of Cauchy problem (\ref{system_u_v}), (\ref{initial_condition_u_v}) where the functions $u(\cdot)$ and $v(\cdot)$ are substituted.
\end{defi}

Note that, according to Definition~\ref{Def_solution}, such a motion is a function $x(\cdot) \in \{x_0\} + I^\alpha (\Linf([0, T], \mathbb{R}^n))$ which, together with $u(\cdot)$ and $v(\cdot),$ satisfies (\ref{system_u_v}) for almost every $t \in [0, T].$

\begin{prop} \label{Prop_motion_u_v_properties}
    For any initial value $x_0 \in \mathbb{R}^n$ and any realizations $u(\cdot) \in \mathcal{U},$ $v(\cdot) \in \mathcal{V},$ there exists the unique motion $x(\cdot) = x(\cdot; x_0, u(\cdot), v(\cdot))$ of system (\ref{system_u_v}), (\ref{initial_condition_u_v}). Moreover, for any $R_0 > 0,$ there exist $\overline{R} > 0$ and $\overline{H} > 0$ such that, for any $x_0 \in B(R_0),$ $u(\cdot) \in \mathcal{U}$ and $v(\cdot) \in \mathcal{V},$ the motion $x(\cdot) = x(\cdot; x_0, u(\cdot), v(\cdot))$ satisfies the following inequalities:
    \begin{equation} \label{Prop_motion_u_v_properties_main}
        \|x(t)\| \leqslant \overline{R}, \quad
        \|x(t) - x(\tau)\| \leqslant \overline{H} |t - \tau|^\alpha, \quad t, \tau \in [0, T].
    \end{equation}
\end{prop}
\proof
    This proposition follows immediately from Theorem~\ref{Thm_existence} and Proposition~\ref{Prop_motion_properties} if we take into account that, for any $u(\cdot) \in \mathcal{U}$ and $v(\cdot) \in \mathcal{V},$ due to ($g.1$)--($g.3$) the function $\overline{f}(t, x) = g(t, x, u(t), v(t)),$ $t \in [0, T],$ $x \in \mathbb{R}^n,$ satisfies the conditions ($f.1$)--($f.3$), and moreover, ($f.3$) is fulfilled with the constant $c_g$ that does not depend on $u(\cdot)$ and $v(\cdot).$
\proofend

Let us consider a guide (see, e.g. \cite[\S~8.2]{NNKrasovskii_AISubbotin_1988}) which is in a certain sense a copy of system (\ref{system_u_v}), (\ref{initial_condition_u_v}). Thus, a motion of the guide is described by the fractional differential equation
\begin{equation} \label{system_y_u_v}
    \begin{array}{c}
        ({}^C D^\alpha y) (t) = g(t, y(t), \widetilde{u}(t), \widetilde{v}(t)), \quad t \in [0, T], \\[0.5em]
        y(t) \in \mathbb{R}^n, \quad \widetilde{u}(t) \in P, \quad \widetilde{v}(t) \in Q,
    \end{array}
\end{equation}
with the initial condition
\begin{equation}\label{initial_condition_y_u_v}
    y(0) = y_0, \quad y_0 \in \mathbb{R}^n.
\end{equation}

Here $y$ is the state vector of the guide, $\widetilde{u}$ and $\widetilde{v}$ are control vectors of the guide; $y_0$ is the initial value. By analogy with Definition~\ref{Def_x} we define a motion $y(\cdot)$ of guide (\ref{system_y_u_v}), (\ref{initial_condition_y_u_v}) that corresponds to an initial value $y_0 \in \mathbb{R}^n$ and admissible realizations $\widetilde{u}(\cdot) \in \mathcal{U},$ $\widetilde{v}(\cdot) \in \mathcal{V}.$ Therefore, from Proposition~\ref{Prop_motion_u_v_properties} it follows that such a motion $y(\cdot) = y(\cdot; y_0, \widetilde{u}(\cdot), \widetilde{v}(\cdot))$ exists and is unique, and moreover, it satisfies estimates similar to (\ref{Prop_motion_u_v_properties_main}).

In the next section, a mutual aiming procedure between initial system (\ref{system_u_v}), (\ref{initial_condition_u_v}) and guide (\ref{system_y_u_v}), (\ref{initial_condition_y_u_v}) is proposed. This procedure is based on the extremal shift rule (see, e.g., \cite[\S\S~2.4, 8.2]{NNKrasovskii_AISubbotin_1988}) and specifies the way of forming control realizations $u(\cdot) \in \mathcal{U}$ and $\widetilde{v}(\cdot) \in \mathcal{V}$ that guarantees proximity between motions of the systems for any disturbance realization $v(\cdot) \in \mathcal{V}$ and any control realization $\widetilde{u}(\cdot) \in \mathcal{U}.$

\section{Mutual Aiming Procedure}

\setcounter{section}{5}
\setcounter{equation}{0}\setcounter{theorem}{0}

Let $x_0, y_0 \in \mathbb{R}^n$ and
\begin{equation*}
    \begin{array}{c}
        \Delta = \{\tau_j\}_{j = 1}^{k + 1} \subset [0, T], \ \
        \tau_1 = 0, \ \  \tau_{j+1} > \tau_j, \ \ j \in \overline{1, k}, \ \ \tau_{k+1} = T, \ \ k \in \mathbb{N},
    \end{array}
\end{equation*}
be a partition of the segment $[0, T].$ Let us consider the following procedure of forming realizations $u(\cdot) \in \mathcal{U}$ in system (\ref{system_u_v}), (\ref{initial_condition_u_v}) and $\widetilde{v}(\cdot) \in \mathcal{V}$ in guide (\ref{system_y_u_v}), (\ref{initial_condition_y_u_v}). Let $j \in \overline{1, k}$ and values $x(\tau_j),$ $y(\tau_j)$ have already been realized. Then we define
\begin{equation} \label{procedure}
    \begin{array}{c}
        \displaystyle
        u(t) = u_j \in \argmin{u \in P} \max_{v \in Q}
        \langle s(\tau_j), g(\tau_j, x(\tau_j), u, v) \rangle, \\[0.5em]
        \displaystyle
        \widetilde{v}(t) = \widetilde{v}_j \in \argmax{\widetilde{v} \in Q} \min_{\widetilde{u} \in P}
        \langle s(\tau_j), g(\tau_j, x(\tau_j), \widetilde{u}, \widetilde{v}) \rangle,
    \end{array} \quad t \in [\tau_j, \tau_{j+1}),
\end{equation}
where we denote
\begin{equation}\label{s}
    s(t) = x(t) - y(t), \quad t \in [0, T].
\end{equation}

\begin{theo} \label{Thm_closeness}
    For any $R_0 > 0$ and $\varepsilon > 0,$ there exist $K > 0$ and $\delta > 0$ such that, for any initial values $x_0, y_0 \in B(R_0),$ any partition $\Delta$ with the diameter $\diam (\Delta) = \max_{j \in \overline{1,k}} (\tau_{j+1} - \tau_j) \leqslant \delta$ and any realizations $v(\cdot) \in \mathcal{V},$ $\widetilde{u}(\cdot) \in \mathcal{U},$ if realizations $u(\cdot) \in \mathcal{U},$ $\widetilde{v}(\cdot) \in \mathcal{V}$ are formed according to mutual aiming procedure (\ref{procedure}), then the corresponding motions $x(\cdot) = x(\cdot; x_0, u(\cdot), v(\cdot))$ of system (\ref{system_u_v}), (\ref{initial_condition_u_v}) and $y(\cdot) = y(\cdot; y_0, \widetilde{u}(\cdot), \widetilde{v}(\cdot))$ of guide (\ref{system_y_u_v}), (\ref{initial_condition_y_u_v}) satisfy the inequality below:
    \begin{equation} \label{Thm_closeness_main}
        \|x(\cdot) - y(\cdot)\|_\infty \leqslant \varepsilon + K \|x_0 - y_0\|.
    \end{equation}
\end{theo}
\proof
    Let $R_0 > 0$ and $\varepsilon > 0.$ By the number $R_0,$ let us choose $\overline{R}$ and $\overline{H}$ in accordance with Proposition~\ref{Prop_motion_u_v_properties}. Due to ($g.2$) let us choose $\lambda_g$ by the number $\overline{R}.$ Let us define $K  = \sqrt{E_\alpha(2 \lambda_g T^\alpha)}.$ Let $\eta > 0$ satisfy the following inequality:
    \begin{equation*}
        \eta \leqslant \Gamma(\alpha + 1) \varepsilon^2/\big(2 T^\alpha E_\alpha(2 \lambda_g T^\alpha) \big).
    \end{equation*}
    By ($g.1$), let us choose $\delta_1 > 0$ such that, for any $t, \tau \in [0, T],$ $x \in B(\overline{R}),$ $u \in P$ and $v \in Q,$ if $|t - \tau| \leqslant \delta_1,$ then
    \begin{equation*} \label{Thm_closeness_delta_1}
        \|g(t, x, u, v) - g(\tau, x, u, v)\| \leqslant \eta / (16 \overline{R}).
    \end{equation*}
    Let $\delta_2 > 0$ be such that
    \begin{equation*} \label{Thm_closeness_delta_2}
        \delta_2^\alpha \leqslant \min\big\{\eta / \big(8 \overline{H} (1 + \overline{R}) c_g \big), \eta / (16 \overline{R} \lambda_g \overline{H}) \big\},
    \end{equation*}
    where $c_g$ is the constant from ($g.3$). Let us define $\delta = \min\{\delta_1, \delta_2\}.$ Let us show that the numbers $K$ and $\delta$ satisfy the statement of the theorem.

    Let $x_0, y_0 \in B(R_0)$ and a partition $\Delta$ has the diameter $\diam(\Delta) \leqslant \delta.$ Let $v(\cdot) \in \mathcal{V},$ $\widetilde{u}(\cdot) \in \mathcal{U}$ and realizations $u(\cdot) \in \mathcal{U},$ $\widetilde{v}(\cdot) \in \mathcal{V}$ be formed according to aiming procedure (\ref{procedure}). Let $x(\cdot) = x(\cdot; x_0, u(\cdot), v(\cdot))$ and $y(\cdot) = y(\cdot; y_0, \widetilde{u}(\cdot), \widetilde{v}(\cdot))$ be the motions of systems (\ref{system_u_v}), (\ref{initial_condition_u_v}) and (\ref{system_y_u_v}), (\ref{initial_condition_y_u_v}).

    Let $s(\cdot)$ be defined by (\ref{s}). Then we have
    \begin{equation*}
        \begin{array}{c}
            s(\cdot) \in \{x_0 - y_0\} + I^\alpha(\Linf([0, T], \mathbb{R}^n)); \\[0.5em]
            \|s(t)\| \leqslant 2 \overline{R}, \quad \|s(t) - s(\tau)\| \leqslant 2 \overline{H} |t - \tau|^\alpha, \quad t, \tau \in [0, T]; \\[0.5em]
            ({}^C D^\alpha s) (t) = g(t, x(t), u(t), v(t)) - g(t, y(t), \widetilde{u}(t), \widetilde{v}(t)) \text{ for a.e. } t \in [0, T].
        \end{array}
    \end{equation*}
    Let us consider the function $\nu(t) = \|s(t)\|^2 - \|x_0 - y_0\|^2,$ $t \in [0, T].$ Since
    \begin{equation*}
        \nu(t) = \|s(t) - (x_0 - y_0)\|^2 + 2 \langle x_0 - y_0, s(t) - (x_0 - y_0) \rangle, \quad t \in [0, T],
    \end{equation*}
    then, due to Corollary~\ref{Cor_differentiating_of_square}, we have $\nu(\cdot) \in I^\alpha(\Linf([0, T], \mathbb{R}))$ and
    \begin{equation} \label{Thm_closeness_estimate}
        \begin{array}{c}
            \hspace*{-0.3em}(D^\alpha \nu)(t)
            \leqslant 2 \langle s(t) - (x_0 - y_0), ({}^C D^\alpha s)(t) \rangle
            + 2 \langle x_0 - y_0, ({}^C D^\alpha s)(t)\rangle \\[0.5em]
            = 2 \langle s(t), ({}^C D^\alpha s)(t) \rangle \text{ for a.e. } t \in [0, T].
        \end{array}
    \end{equation}

    Let us show that
    \begin{equation} \label{Thm_closeness_ineq_base}
        \langle s(t), ({}^C D^\alpha s)(t) \rangle \leqslant \lambda_g \|s(t)\|^2 + \eta \text{ for a.e. } t \in [0, T].
    \end{equation}
    For almost every $t \in [0, T],$ we obtain
    \begin{equation} \label{Thm_closeness_ineq_1}
        \begin{array}{c}
            \hspace*{-0.3em}\langle s(t), ({}^C D^\alpha s)(t) \rangle
            = \langle s(t), g(t, x(t), u(t), v(t)) - g(t, x(t), \widetilde{u}(t), \widetilde{v}(t)) \rangle \\[0.5em]
            + \langle s(t), g(t, x(t), \widetilde{u}(t), \widetilde{v}(t)) - g(t, y(t), \widetilde{u}(t), \widetilde{v}(t)) \rangle.
        \end{array}
    \end{equation}
    Let us estimate each of the two terms separately.

    Let $j \in \overline{1, k}$ and $t \in [\tau_j, \tau_{j+1}).$ By $(g.3)$ and the choice of $\delta_2,$ we derive
    \begin{equation*}
        \begin{array}{c}
            \langle s(t), g(t, x(t), u(t), v(t)) - g(t, x(t), \widetilde{u}(t), \widetilde{v}(t)) \rangle \\[0.5em]
            \leqslant \langle s(\tau_j), g(t, x(t), u(t), v(t)) - g(t, x(t), \widetilde{u}(t), \widetilde{v}(t)) \rangle \\[0.5em]
            + \|s(t) - s(\tau_j)\| \big(\|g(t, x(t), u(t), v(t))\| + \|g(t, x(t), \widetilde{u}(t), \widetilde{v}(t))\| \big) \\[0.5em]
            \leqslant \langle s(\tau_j), g(t, x(t), u(t), v(t)) - g(t, x(t), \widetilde{u}(t), \widetilde{v}(t)) \rangle 
            + 4 \overline{H} (1 + \overline{R}) c_g \delta^\alpha \\[0.5em]
            \leqslant \langle s(\tau_j), g(t, x(t), u(t), v(t)) - g(t, x(t), \widetilde{u}(t), \widetilde{v}(t)) \rangle + \eta/2.
        \end{array}
    \end{equation*}
    Further, due to the choice of $\lambda_g,$ $\delta_1$ and $\delta_2$ we obtain
    \begin{equation*}
        \begin{array}{c}
            \langle s(\tau_j), g(t, x(t), u(t), v(t)) \rangle 
            \leqslant \langle s(\tau_j), g(\tau_j, x(\tau_j), u(t), v(t)) \rangle \\[0.5em]
            + \|s(\tau_j)\| \|g(t, x(t), u(t), v(t)) - g(\tau_j, x(t), u(t), v(t)) \| \\[0.5em]
            + \|s(\tau_j)\| \|g(\tau_j, x(t), u(t), v(t)) - g(\tau_j, x(\tau_j), u(t), v(t)) \| \\[0.5em]
            \leqslant \langle s(\tau_j), g(\tau_j, x(\tau_j), u(t), v(t)) \rangle \\[0.5em]
            + 2 \overline{R} \|g(t, x(t), u(t), v(t)) - g(\tau_j, x(t), u(t), v(t)) \| 
            + 2 \overline{R} \lambda_g \overline{H} \delta^\alpha \\[0.5em]
            \leqslant \langle s(\tau_j), g(\tau_j, x(\tau_j), u(t), v(t)) \rangle + \eta/4,
        \end{array}
    \end{equation*}
    and, similarly,
    \begin{equation*}
        \begin{array}{c}
            \langle s(\tau_j), g(t, x(t), \widetilde{u}(t), \widetilde{v}(t))
            \geqslant \langle s(\tau_j), g(\tau_j, x(\tau_j), \widetilde{u}(t), \widetilde{v}(t)) \rangle - \eta/4.
        \end{array}
    \end{equation*}
    Finally, in accordance with ($g.4$) and choice (\ref{procedure}) of $u_j,$ $\widetilde{v}_j$ we get
    \begin{equation*}
        \begin{array}{c}
            \langle s(\tau_j), g(\tau_j, x(\tau_j), u(t), v(t)) \rangle
            - \langle s(\tau_j), g(\tau_j, x(\tau_j), \widetilde{u}(t), \widetilde{v}(t)) \rangle \\[0.5em]
            = \langle s(\tau_j), g(\tau_j, x(\tau_j), u_j, v(t)) \rangle
            - \langle s(\tau_j), g(\tau_j, x(\tau_j), \widetilde{u}(t), \widetilde{v}_j) \rangle \\[0.5em]
            \leqslant \max\limits_{v \in Q} \langle s(\tau_j), g(\tau_j, x(\tau_j), u_j, v) \rangle
            - \min\limits_{\widetilde{u} \in P} \langle s(\tau_j), g(\tau_j, x(\tau_j), \widetilde{u}, \widetilde{v}_j) \rangle \\[0.5em]
            = \min\limits_{u \in P} \max\limits_{v \in Q}
            \langle s(\tau_j), g(\tau_j, x(\tau_j), u, v) \rangle
            - \max\limits_{\widetilde{v} \in Q} \min\limits_{\widetilde{u} \in P}
            \langle s(\tau_j), g(\tau_j, x(\tau_j), \widetilde{u}, \widetilde{v}) \rangle
            = 0.
        \end{array}
    \end{equation*}
    Consequently, for $t \in [0, T),$ we have
    \begin{equation}\label{Thm_closeness_ineq_1.5}
        \langle s(t), g(t, x(t), u(t), v(t)) - g(t, x(t), \widetilde{u}(t), \widetilde{v}(t)) \rangle
        \leqslant \eta.
    \end{equation}

    Let us estimate the second term in (\ref{Thm_closeness_ineq_1}). For $t \in [0, T],$ due to the choice of $\lambda_g$ we derive
    \begin{equation} \label{Thm_closeness_ineq_2}
        \begin{array}{c}
            \langle s(t), g(t, x(t), \widetilde{u}(t), \widetilde{v}(t))
            - g(t, y (t), \widetilde{u}(t), \widetilde{v}(t)) \rangle \\[0.5em]
            \leqslant \|s(t)\| \| g(t, x(t), \widetilde{u}(t), \widetilde{v}(t))
            - g(t, y(t), \widetilde{u}(t), \widetilde{v}(t)) \|
            \leqslant \lambda_g \|s(t)\|^2.
        \end{array}
    \end{equation}

    Thus, the validity of inequality (\ref{Thm_closeness_ineq_base}) follows from (\ref{Thm_closeness_ineq_1})--(\ref{Thm_closeness_ineq_2}).

    From (\ref{Thm_closeness_estimate}) and (\ref{Thm_closeness_ineq_base}) we obtain
    \begin{equation*}
        (D^\alpha \nu)(t)
        \leqslant 2 \lambda_g \|s(t)\|^2 + 2 \eta \text{ for a.e. } t \in [0, T].
    \end{equation*}
    Therefore, according to ($B.2$), for every $t \in [0, T],$ we have
    \begin{equation*}
        \nu(t) \leqslant \frac{1}{\Gamma(\alpha)} \int_{0}^{t} \frac{2 \lambda_g \|s(\tau)\|^2 + 2 \eta}{(t - \tau)^{1 - \alpha}} d \tau
        \leqslant \frac{2 \eta T^\alpha}{\Gamma(\alpha + 1)}
        + \frac{2 \lambda_g}{\Gamma(\alpha)} \int_{0}^{t} \frac{\|s(\tau)\|^2}{(t - \tau)^{1 - \alpha}} d \tau.
    \end{equation*}
    Consequently, due to the definition of $\nu(\cdot)$ we deduce
    \begin{equation*}
        \|s(t)\|^2
        \leqslant \frac{2 \eta T^\alpha}{\Gamma(\alpha + 1)}
        + \|x_0 - y_0\|^2
        + \frac{2 \lambda_g}{\Gamma(\alpha)} \int_{0}^{t} \frac{\|s(\tau)\|^2}{(t - \tau)^{1 - \alpha}} d \tau.
    \end{equation*}
    Hence, by Lemma~\ref{Lem_Bellman_Gronwall} and the choice of $\eta$ and $K,$ for $t \in [0, T],$ we obtain
    \begin{equation*}
        \|s(t)\|^2 \leqslant \Big( \frac{2 \eta T^\alpha}{\Gamma(\alpha + 1)} + \|x_0 - y_0\|^2 \Big) E_\alpha(2 \lambda_g T^\alpha)
        \leqslant \varepsilon^2 + K^2 \|x_0 - y_0\|^2.
    \end{equation*}
    Thus, inequality (\ref{Thm_closeness_main}) and the theorem are proved.
\proofend

\section{Example}

\setcounter{section}{6}
\setcounter{equation}{0}\setcounter{theorem}{0}

Let us illustrate the constructions from Sect.~4 and 5 by an example. Let a motion of the conflict-controlled dynamical system be described by the fractional differential equations
\begin{equation} \label{system_u_v_ex}
    \begin{array}{c}
        \begin{cases}
            (^C D^{0.5} x_1) (t) = x_2(t) + 0.3 u_1(t) + 0.4 v_1(t), \\[0.5em]
            (^C D^{0.5} x_2) (t) = - \sin(x_1(t)) + \cos(t) + 0.5 u_2(t) + 0.2 v_2(t),
        \end{cases} \\[1.5em]
        t \in [0, 5], \quad x(t) = (x_1(t), x_2(t)) \in \mathbb{R}^2, \\[0.5em]
        u(t) = (u_1(t), u_2(t)) \in P = \{u \in \mathbb{R}^2: \|u\| \leqslant 1\}, \\[0.5em]
        v(t) = (v_1(t), v_2(t)) \in Q = \{v \in \mathbb{R}^2: \|v\| \leqslant 1\},
    \end{array}
\end{equation}
with the initial condition
\begin{equation}\label{initial_condition_u_v_ex}
    x(0) = (-1, 0).
\end{equation}

Let us consider a guide which motion is described by the similar fractional differential equations
\begin{equation} \label{system_y_u_v_ex}
    \begin{array}{c}
        \begin{cases}
            (^C D^{0.5} y_1) (t) = y_2(t) + 0.3 \widetilde{u}_1(t) + 0.4 \widetilde{v}_1(t), \\[0.5em]
            (^C D^{0.5} y_2) (t) = - \sin(y_1(t)) + \cos(t) + 0.5 \widetilde{u}_2(t) + 0.2 \widetilde{v}_2(t),
        \end{cases} \\[1.5em]
        t \in [0, 5], \quad y(t) = (y_1(t), y_2(t)) \in \mathbb{R}^2, \\[0.5em]
        \widetilde{u}(t) = (\widetilde{u}_1(t), \widetilde{u}_2(t)) \in P, \quad
        \widetilde{v}(t) = (\widetilde{v}_1(t), \widetilde{v}_2(t)) \in Q,
    \end{array}
\end{equation}
with the initial condition
\begin{equation}\label{initial_condition_y_u_v_ex}
    y(0) = (0, 1).
\end{equation}

\begin{center}

\begingroup
  \makeatletter
  \providecommand\color[2][]{%
    \GenericError{(gnuplot) \space\space\space\@spaces}{%
      Package color not loaded in conjunction with
      terminal option `colourtext'%
    }{See the gnuplot documentation for explanation.%
    }{Either use 'blacktext' in gnuplot or load the package
      color.sty in LaTeX.}%
    \renewcommand\color[2][]{}%
  }%
  \providecommand\includegraphics[2][]{%
    \GenericError{(gnuplot) \space\space\space\@spaces}{%
      Package graphicx or graphics not loaded%
    }{See the gnuplot documentation for explanation.%
    }{The gnuplot epslatex terminal needs graphicx.sty or graphics.sty.}%
    \renewcommand\includegraphics[2][]{}%
  }%
  \providecommand\rotatebox[2]{#2}%
  \@ifundefined{ifGPcolor}{%
    \newif\ifGPcolor
    \GPcolorfalse
  }{}%
  \@ifundefined{ifGPblacktext}{%
    \newif\ifGPblacktext
    \GPblacktexttrue
  }{}%
  \let\gplgaddtomacro\g@addto@macro
  \gdef\gplbacktext{}%
  \gdef\gplfronttext{}%
  \makeatother
  \ifGPblacktext
    \def\colorrgb#1{}%
    \def\colorgray#1{}%
  \else
    \ifGPcolor
      \def\colorrgb#1{\color[rgb]{#1}}%
      \def\colorgray#1{\color[gray]{#1}}%
      \expandafter\def\csname LTw\endcsname{\color{white}}%
      \expandafter\def\csname LTb\endcsname{\color{black}}%
      \expandafter\def\csname LTa\endcsname{\color{black}}%
      \expandafter\def\csname LT0\endcsname{\color[rgb]{1,0,0}}%
      \expandafter\def\csname LT1\endcsname{\color[rgb]{0,1,0}}%
      \expandafter\def\csname LT2\endcsname{\color[rgb]{0,0,1}}%
      \expandafter\def\csname LT3\endcsname{\color[rgb]{1,0,1}}%
      \expandafter\def\csname LT4\endcsname{\color[rgb]{0,1,1}}%
      \expandafter\def\csname LT5\endcsname{\color[rgb]{1,1,0}}%
      \expandafter\def\csname LT6\endcsname{\color[rgb]{0,0,0}}%
      \expandafter\def\csname LT7\endcsname{\color[rgb]{1,0.3,0}}%
      \expandafter\def\csname LT8\endcsname{\color[rgb]{0.5,0.5,0.5}}%
    \else
      \def\colorrgb#1{\color{black}}%
      \def\colorgray#1{\color[gray]{#1}}%
      \expandafter\def\csname LTw\endcsname{\color{white}}%
      \expandafter\def\csname LTb\endcsname{\color{black}}%
      \expandafter\def\csname LTa\endcsname{\color{black}}%
      \expandafter\def\csname LT0\endcsname{\color{black}}%
      \expandafter\def\csname LT1\endcsname{\color{black}}%
      \expandafter\def\csname LT2\endcsname{\color{black}}%
      \expandafter\def\csname LT3\endcsname{\color{black}}%
      \expandafter\def\csname LT4\endcsname{\color{black}}%
      \expandafter\def\csname LT5\endcsname{\color{black}}%
      \expandafter\def\csname LT6\endcsname{\color{black}}%
      \expandafter\def\csname LT7\endcsname{\color{black}}%
      \expandafter\def\csname LT8\endcsname{\color{black}}%
    \fi
  \fi
    \setlength{\unitlength}{0.0500bp}%
    \ifx\gptboxheight\undefined%
      \newlength{\gptboxheight}%
      \newlength{\gptboxwidth}%
      \newsavebox{\gptboxtext}%
    \fi%
    \setlength{\fboxrule}{0.5pt}%
    \setlength{\fboxsep}{1pt}%
\begin{picture}(5760.00,3456.00)%
    \gplgaddtomacro\gplbacktext{%
      \csname LTb\endcsname%
      \put(726,464){\makebox(0,0)[r]{\strut{}$-1$}}%
      \csname LTb\endcsname%
      \put(726,1365){\makebox(0,0)[r]{\strut{}$0$}}%
      \csname LTb\endcsname%
      \put(726,2266){\makebox(0,0)[r]{\strut{}$1$}}%
      \csname LTb\endcsname%
      \put(726,3167){\makebox(0,0)[r]{\strut{}$x, y$}}%
      \csname LTb\endcsname%
      \put(858,244){\makebox(0,0){\strut{}$0$}}%
      \csname LTb\endcsname%
      \put(1759,244){\makebox(0,0){\strut{}$1$}}%
      \csname LTb\endcsname%
      \put(2660,244){\makebox(0,0){\strut{}$2$}}%
      \csname LTb\endcsname%
      \put(3561,244){\makebox(0,0){\strut{}$3$}}%
      \csname LTb\endcsname%
      \put(4462,244){\makebox(0,0){\strut{}$4$}}%
      \csname LTb\endcsname%
      \put(5363,244){\makebox(0,0){\strut{}$t$}}%
    }%
    \gplgaddtomacro\gplfronttext{%
      \csname LTb\endcsname%
      \put(5183,2900){\makebox(0,0)[r]{\strut{}$x_1$}}%
      \csname LTb\endcsname%
      \put(5183,2636){\makebox(0,0)[r]{\strut{}$x_2$}}%
      \csname LTb\endcsname%
      \put(5183,2372){\makebox(0,0)[r]{\strut{}$y_1$}}%
      \csname LTb\endcsname%
      \put(5183,2108){\makebox(0,0)[r]{\strut{}$y_2$}}%
    }%
    \gplbacktext
    \put(0,0){\includegraphics{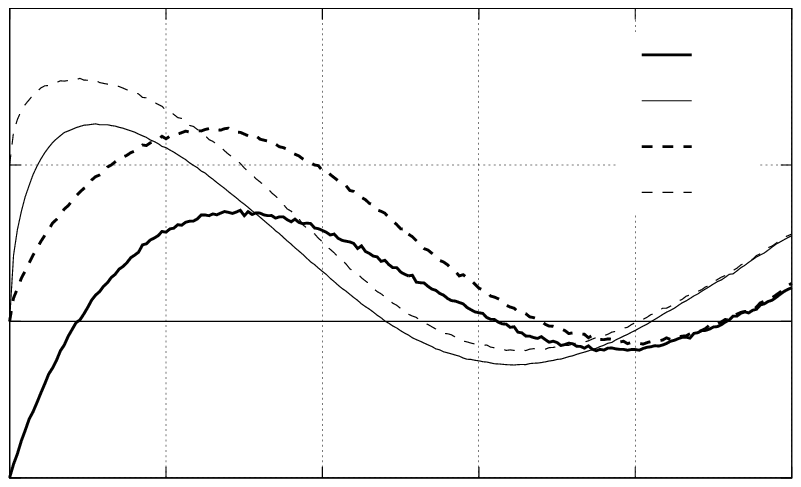}}%
    \gplfronttext
  \end{picture}%
\endgroup

    \bigskip

    Fig. 6.1: The realized motions of system (\ref{system_u_v_ex}), (\ref{initial_condition_u_v_ex}) and guide (\ref{system_y_u_v_ex}), (\ref{initial_condition_y_u_v_ex}).
\end{center}

For system (\ref{system_u_v_ex}), (\ref{initial_condition_u_v_ex}) and guide (\ref{system_y_u_v_ex}), (\ref{initial_condition_y_u_v_ex}), mutual aiming procedure (\ref{procedure}) was simulated. The uniform partition $\Delta$ of the segment $[0, 5]$ with the step $\delta = 0.0005$ was chosen. Realizations $u(\cdot)$ in the initial system and $\widetilde{v}(\cdot)$ in the guide were formed according to procedure (\ref{procedure}), while realizations $v(\cdot)$ in the initial system and $\widetilde{u}(\cdot)$ in the guide were formed as piecewise constant on the partition $\Delta$ functions with random values from $P$ and $Q,$ respectively. For the numerical simulation of motions of system (\ref{system_u_v_ex}), (\ref{initial_condition_u_v_ex}) and guide (\ref{system_y_u_v_ex}), (\ref{initial_condition_y_u_v_ex}) the fractional forward Euler method (see, e.g., \cite[p.~101]{CLi_FZeng_2015}) was used. The obtained results, presented in Fig.~6.1, show that the realized motions $x(\cdot)$ and $y(\cdot)$ of the systems are close to each other despite the choice of the realizations $v(\cdot)$ and $\widetilde{u}(\cdot),$ which agrees with Theorem~\ref{Thm_closeness}.

\section{Conclusion}

\setcounter{section}{7}
\setcounter{equation}{0}\setcounter{theorem}{0}

In the paper a conflict-controlled dynamical system described by an ordinary fractional differential equation with the Caputo derivative of an order $\alpha \in (0, 1)$ is considered. A suitable notion of a system motion that does not assume its differentiability is proposed. The existence and uniqueness results for such a motion are obtained. An auxiliary guide is introduced which is in a certain sense a copy of the initial system. In order to ensure proximity between motions of the system and the guide a mutual aiming procedure is elaborated. To justify this aiming procedure, the estimate of the fractional derivative of the superposition of a convex Lyapunov function and a motion of the system is proved. The obtained results are illustrated by an example.

Let us stress again that the proposed aiming procedure guarantees proximity between initial system (\ref{system_u_v}), (\ref{initial_condition_u_v}) and guide (\ref{system_y_u_v}), (\ref{initial_condition_y_u_v}) for any disturbances $v(t)$ and any control actions $\widetilde{u}(t).$ Therefore, in the further applications control actions $\widetilde{u}(t)$ in the guide may be used in order to compensate disturbances $v(t)$ and ensure the desired quality of a control process in the initial system.

\section*{Acknowledgements}

This work is supported by the Program of the Presidium of the Russian Academy of Sciences No. 01 'Fundamental Mathematics and its Applications' under grant PRAS-18-01.




 \bigskip \smallskip

 \it

 \noindent
$^1$ Krasovskii Institute of Mathematics and Mechanics \\
Ural Branch of Russian Academy of Sciences \\
S. Kovalevskaya Str., Block 16 \\
Ekaterinburg -- 620990, RUSSIA  \\[4pt]
e-mail: m.i.gomoyunov@gmail.com
\hfill Received: October 5, 2017 \\[12pt]

\end{document}